\documentclass[12pt]{amsart}

\usepackage{xy, graphicx, color,hyperref,verbatim}
\pdfoutput=1
\usepackage{amsfonts, ytableau,array}
\usepackage[utf8]{inputenc}
\usepackage{amssymb,amsmath,amsthm,bm,mathrsfs,enumerate,bbm,enumitem}
\usepackage[a4paper, total={6in, 8in}]{geometry}
\usepackage[english]{babel}
\usepackage{fancyhdr}
\usepackage[utf8]{inputenc}
\usepackage{mathtools,cite,mathdots}
\usepackage{tikz}
\usepackage{tikz-cd}
\usetikzlibrary{matrix}
\usepackage{siunitx}
\usepackage{booktabs}
\usepackage{blindtext}
\usepackage{array}
\newcolumntype{P}[1]{>{\centering\arraybackslash}p{#1}}
\PassOptionsToPackage{linktocpage}{hyperref}
   \hypersetup{
           breaklinks=true,  
           colorlinks=true,  
           pdfusetitle=true, 
           citecolor={blue!80!black}
        }
\makeatletter
\def\subsection{\@startsection{subsection}{3}
  \z@{.9\linespacing\@plus.7\linespacing}{.1\linespacing}
  {\normalfont\bfseries}}
\makeatother

 \newtheorem{thm}{Theorem}
\newtheorem{lemma}{Lemma}

\newtheorem{prop}{Proposition}

\newtheorem{cor}{Corollary}[thm]
 
\theoremstyle{definition}
\newtheorem{remark}{Remark}

\newcommand{\nc}{\newcommand}

\nc{\mc}{\mathcal}
\nc{\mb}{\mathbb}
\nc{\mf}{\mathfrak}
\nc{\ul}{\underline}
\nc{\ol}{\overline}
\nc{\N}{\mb N}
\nc{\R}{\mb R}
\nc{\Z}{\mb Z}
\nc{\Q}{\mb Q}
\nc{\C}{\mb C}
\nc{\F}{\mb F}
\nc{\dmo}{\DeclareMathOperator}

\nc{\mat}[4]{
    \begin{pmatrix}
      #1 & #2 \\
      #3 & #4
    \end{pmatrix}
}
\dmo{\Ker}{Ker} \dmo{\val}{val} \dmo{\ord}{ord}
\dmo{\I}{I}
\dmo{\II}{II}
\dmo{\odd}{odd}
\dmo{\sgn}{sgn}
\dmo{\dett}{det}
\dmo{\Span}{Span}

\dmo{\diag}{diag}

\nc{\beq}{\begin{equation*}}
\nc{\eeq}{\end{equation*}}
\nc{\half}{\frac{1}{2}}

\dmo{\Mod}{mod}

\dmo{\core}{core}
\dmo{\res}{res}
\dmo{\lin}{lin}
 
 \dmo{\St}{St}
 \dmo{\st}{st}
 
 \dmo{\Tr}{Tr}
 \dmo{\RO}{RO}
\dmo{\Sp}{Sp}
\dmo{\SO}{SO}
\dmo{\SL}{SL}
\dmo{\GL}{GL}
 \dmo{\Spin}{Spin}
\dmo{\GSp}{GSp}
\nc{\la}{\lambda}
  \nc{\eps}{\varepsilon}
  
 \nc{\lip}{\langle}
 \nc{\rip}{\rangle}
\nc{\gm}{\gamma}

\dmo{\Perm}{Perm}
\dmo{\Res}{Res}
\dmo{\Ind}{Ind}
\dmo{\ind}{ind}
\dmo{\tr}{tr}
\dmo{\Sym}{Sym}
\dmo{\reg}{reg}
\dmo{\End}{End}
\dmo{\Hom}{Hom}
\dmo{\Pin}{Pin}
\dmo{\Or}{O}
\dmo{\SW}{SW}
\title[Stiefel-Whitney Classes]{Stiefel-Whitney Classes of representations of some finite groups of lie type} 
\author{Neha Malik}
\author{Steven Spallone}

\DeclareMathOperator{\pr}{pr}

\dmo{\ORep}{ORep}
\dmo{\Rep}{Rep}

\usepackage{blindtext,graphicx}
\usepackage[absolute]{textpos}
\newcommand{\f}{\mathbb{F}}

\newcommand{\rr}{\mathbb{R}}

\address{Indian Institute of Science Education and Research, Pune-411008,Maharashtra,India}
\email{nehamalik5194@gmail.com}
\address{Indian Institute of Science Education and Research, Pune-411008,Maharashtra,India}
\email{sspallone@gmail.com}
\date{\today}

\begin{document}

\begin{abstract}
  In this note we present the Stiefel-Whitney classes (SWCs) for orthogonal representations of several finite groups of Lie type, namely for $G=\SL(2,q),$ $\SL(3,q),$ $\Sp(4,q)$, and $\Sp(6,q)$, with $q$ odd. We also describe the SWCs for $G=\SL(2,q)$ when $q$ is even.
\end{abstract}
\maketitle
\tableofcontents
\section{Introduction}

The theory of Stiefel-Whitney classes (SWCs) of vector bundles is an old unifying concept in geometry.  Real representations, or equivalently, orthogonal complex representations of a group $G$ give rise to \emph{flat} vector bundles over the classifying space $BG$. Via this construction one defines SWCs of representations, which live in the group cohomology $H^*(G,\Z/2\Z)$.  There do not seem to be many explicit calculations of such SWCs, especially when $G$ is nonabelian. In this note we give formulas for SWCs in terms of character values, for the groups mentioned in the abstract. 

Formulas for cyclic groups are well-known; we review this below. The second SWCs of representations of $S_n$ and related groups were found in \cite{Ganguly}.
The case of $\GL(2,q)$ with $q$ odd was completed recently in \cite{GJ}; their results are analogous to ours.

This note is intended only as an announcement of our calculations; the work will appear elsewhere.\\

\textbf{Acknowledgements.} 
The authors would like to thank Ryan Vinroot, Dipendra Prasad, and Rohit Joshi for helpful discussions.  The research by the first author is supported by a PhD fellowship from the Council of Scientific \& Industrial Research, India.

\section{Notations \& Preliminaries}
\subsection{Orthogonal Representations}

Let $G$ be a finite group. We say that a representation $(\pi,V)$, with $V$ a complex finite-dimensional vector space, is \emph{orthogonal}, provided there exists a  non-degenerate $G$-invariant symmetric bilinear form $B:V\times V\to \C$. From \cite[Chapter II, Section 6]{BrokerDieck}, we note:

\begin{prop} \label{prop.one}
If $(\pi,V)$ is orthogonal, then there is a $(\pi_0,V_0)$, with $V_0$ a real vector space, so that $\pi_0 \otimes_{\R} \C \cong \pi$. Moreover, $\pi_0$ is unique up to isomorphism.
\end{prop}

\subsection{Stiefel-Whitney Classes}

In this section we define Stiefel-Whitney classes for orthogonal complex representations of finite groups. This is equivalent to the definition in \cite{Benson} or \cite{gunarwardena1989stiefel}, given for real representations.

Recall from  \cite{milnor} that a ($d$-dimensional real) vector bundle $\mc V$ over a paracompact space $B$ is assigned a sequence 
\beq
w_1(\mc V), \ldots, w_d(\mc V)
\eeq
of cohomology classes, with each $w_i(\mc V) \in H^*(B,\Z/2\Z)$.

For every finite group $G$ there is a classifying space $BG$ with a principal $G$-bundle $EG$, unique up to homotopy. From a representation $(\pi,V)$ one can form the associated vector bundle $EG[V]$ over $BG$. Then one puts
\beq
w_i^\R(\pi)=w_i(EG[V]);
\eeq
see for instance \cite{Benson} or  \cite{gunarwardena1989stiefel}. 
Moreover the singular cohomology $H^*(BG,\Z/2\Z)$ is isomorphic to the group cohomology $H^*(G,\Z/2\Z)$. 

Hence, from a real representation $\rho$, we can associate cohomology classes $w_i^\R(\rho) \in H^i(G,\Z/2\Z)$, called Stiefel-Whitney classes (SWCs).
We prefer to work with an orthogonal complex representation $\pi$, which by Proposition  \ref{prop.one} comes from a unique real representation $\pi_0$. 
We thus define
$$w_i(\pi):=w_i^\rr(\pi_0),$$
for $0\leq i\leq \deg\pi$.
 
 At this point we simply write $H^*(G)=H^*(G,\Z/2\Z)$.  A group homomorphism $\varphi: G_1 \to G_2$ induces a map 
  \beq
  \varphi^*: H^*(G_2) \to H^*(G_1).
  \eeq

 The first few $w_i(\pi)$ have nice interpretations. Firstly, $w_0(\pi)=1 \in H^0(G)$. Moreover, the first SWC, applied to linear characters $G \to \{\pm 1\} \cong \Z/2\Z$, is the well-known isomorphism
  \beq
  w_1: \Hom(G, \Z/2\Z)  \overset{\sim}{\to} H^1(G,\Z/2\Z).
  \eeq

More generally, if $\pi$ is an orthogonal representation, then $w_1(\pi)=w_1(\det \pi)$, where $\det \pi$ is simply the composition of $\pi$ with the determinant map.
When $\det \pi=1$, then $w_2(\pi)$ vanishes iff $\pi$ lifts to the corresponding spin group. (See \cite{Ganguly} for details.)
  
  The SWCs are functorial in the following sense: If  $\rho$ is an orthogonal representation of $G_2$, then  
  \begin{equation} \label{funct.w}
  \varphi^*(w(\pi))=w(\pi \circ \varphi).
  \end{equation}
  
\section{Group Cohomology}  
 
 In this section we review explicit descriptions of group cohomology for several important groups.   Later, we will describe SWCs in terms of these descriptions.
  Throughout this note, we put $k=\Z/2\Z$, and write $\Theta_\pi$ for the character of a representation $\pi$.

\subsection{Cyclic Groups}\label{cy}

Let $n$ be even, and $G=C_n$, the cyclic group of order $n$. Let $g$ be a generator of $G$, and $\chi_\bullet$ the linear character of $G$ with $\chi_\bullet(g)=e^{\frac{2\pi i}{n} }$.   Write `sgn' for the linear character of order 2.
It is known \cite{kamber1967flat} that 
$$H^*(C_n)=\begin{cases}
k[s,t]/(s^2),& n\equiv0 \text{ (mod } 4)\\
k[v],& n\equiv2 \text{ (mod }4)
\end{cases}$$
where $s=w_1(\sgn),$ $t=w_2(\chi_\bullet \oplus \chi_\bullet^{-1})$ for $n\equiv0 \text{ (mod } 4)$, and $v=w_1(\sgn)$ when $n \equiv2 \text{ (mod }4)$. \\

Since the cyclic case is so easy, we can give our first example of expressing SWCs in terms of character values.

\begin{prop} Let $\pi$ be an orthogonal representation of $C_n$, and $b_\pi=\frac{1}{4} \left( \deg \pi-\Theta_\pi(g^{n/2}) \right)$.
\begin{enumerate}
\item If $n \equiv 2 \mod 4$, then
\beq
w(\pi)=(1+v)^{2b_\pi}.
\eeq
\item If $n \equiv 0 \mod 4$, and $\det \pi=1$, then
\beq
w(\pi)=(1+t)^{b_\pi}.
\eeq
\item If $n \equiv 0 \mod 4$, and $\det \pi \neq 1$, then
\beq
w(\pi)=(1+s)(1+t)^{b_\pi}.
\eeq
\end{enumerate}
\end{prop}

Let $C_n^r$ be the $n$-fold product of $C_n$, with projection maps $\pr_i: C_n^r\to C_n$ for $1 \leq i \leq r$. By K\"{u}nneth, we have 
$$H^*(C_n^r)=\begin{cases}
k[s_1,\dots,s_r,t_1,\dots,t_r]/(s_1^2,\hdots,s_r^2),& n\equiv0 \text{ (mod } 4)\\
k[v_1,\dots,v_r],& n\equiv2 \text{ (mod }4)
\end{cases}$$
where   we put $s_i=w_1(\sgn\circ\pr_i)$ and $t_i=w_2((\chi_\bullet \oplus\chi_\bullet^{-1})\circ\pr_i)$  for $n\equiv0 \text{ (mod } 4)$, and $v_i=w_1(\sgn\circ\pr_i)$, for $n\equiv2 \text{ (mod } 4)$.  (Here $1\leq i\leq r$.)

\subsection{Special Linear Groups} \label{sl.section}
%Let $G=\SL(2,q)$ with $q$ odd.
From \cite{Quillen}, we have for odd $q$,
$$H^*(\SL(2,q))=k[e]\otimes_k k[b]/b^2,$$
where $\deg(b)=3$ and  $e$ is the unique non-zero element in $H^4(\SL(2,q))$.\\
Let $\mathcal{Z}$ be the center of $\SL(2,q)$, which has order $2$. The restriction map \begin{align*}
H^4(\SL(2,q))&\to H^4(\mathcal{Z})\\
e&\mapsto v^4
\end{align*}
is a group isomorphism.
\subsection{Symplectic Groups}
 Fix the matrix $J=\begin{pmatrix}
0&&0&&\hdots&&0&&1\\
0&&0&&\hdots&&-1&&0\\
\vdots&&\vdots&&\iddots&&\vdots&&\vdots\\
0&&1&&\hdots&&0&&0\\
-1&&0&&\hdots&&0&&0\\
\end{pmatrix}$,
\\
\vspace{2mm}

and let $G=\Sp(2n,q):=\{A\in \GL(2n,q)\mid A^tJA=J\}$, with $q$ odd. Write $X$ for the subgroup of matrices in $G$, all of whose nonzero entries lie either on the diagonal or the antidiagonal. 
Thus
\beq
X=\begin{pmatrix}
\square&0&0&\hdots&0&0&\square\\
0&\ddots&& &&\iddots&0\\
0&&\square&&\square&&0\\
\vdots&&\vdots&\scalebox{3}{$\square$}_{2\times 2}&\vdots&&\vdots\\
0&&\square&&\square&&0\\
0&\iddots&&&&\ddots&0\\
\square&0&0&\hdots&0&0&\square\\
\end{pmatrix}.
\eeq
Note that $X$ is isomorphic to the direct product of $n$ copies of $\SL(2,q)$.
 So there are projections
$$\text{pr}_j:X\to \SL(2,q) \text{ ; }1\leq j\leq n,$$ and by  K\"{u}nneth we have
$$H^*(X)\cong k[e_1,\hdots,e_n]\otimes_k k[b_1,\hdots,b_n]/(b_1^2, \ldots, b_n^2),$$

where $e_j=\text{pr}^*_j(e)$ and $b_j=\text{pr}^*_j(b)$ for $1\leq j\leq n.$

\section{Detection of Cohomology}
Let $G$ be a finite group. We say a subgroup $G'$ \emph{detects} the mod-2 cohomology of $G$ when the restriction map
$$H^*(G)\to H^*(G')$$
is an injection.

Suppose this is the case, and that $\pi$ is an orthogonal representation of $G$. Then $w(\pi)$ is specified by its image in $H^*(G')$, which is actually the SWC of the restriction of $\pi$ to $G'$. In our formulas in the next section, when there is a detecting subgroup, we will simply again write $w(\pi)$ for this image.

\begin{lemma}[\cite{Milgram}, Cor. 5.2, Ch. II]
A Sylow 2-subgroup of $G$ detects its mod-2 cohomology.  \end{lemma}

\begin{prop}[\cite{Quillen}]
When $q$ is odd, the diagonal subgroup $A$ of $\GL(n,q)$ detects its mod-2 cohomology.
\end{prop}

From the above one can deduce:
\begin{prop}
When $n$ and $q$ are odd, the diagonal subgroup $A_1$ of $\SL(n,q)$ detects its mod-2 cohomology.
\end{prop}

Finally, there is:

\begin{prop} [\cite{Milgram}, Lemma 6.2]
For odd $q$, the subgroup $X$ detects the mod-2 cohomology of $\Sp(2n,q)$. 
\end{prop}

\section{Main Theorems}
 Now we are ready to  describe the promised SWCs.

\subsection{$G=\SL(2,q)$, $q$ odd}

\begin{thm}\label{thm2.section} 
The total SWC of an orthogonal representation $\pi$ of $\SL(2,q)$ is
\begin{equation}
w(\pi)=(1+e)^{r_\pi},
\end{equation}
where $e$ is as defined in Section \ref{sl.section}, and $r_\pi=\frac{1}{8}(\Theta_\pi(\mathbbm{1})-\Theta_\pi(-\mathbbm{1})).$
\end{thm}

We can say more when $\pi$ is irreducible; write $\omega_\pi$ for its central character.

\begin{prop}[\cite{gow1985real}, Theorem 1]\label{gow} If $\pi$ is irreducible orthogonal, then $\omega_\pi(- \mathbbm{1})=1$. If $\pi$ is irreducible symplectic, then $\omega_\pi(- \mathbbm{1})=-1$.
\end{prop}

Thus if $\pi$ is irreducible orthogonal, then $\Theta_\pi(-\mathbbm{1})=\Theta_\pi(\mathbbm{1})$, and so $r_\pi=0$. Therefore:

\begin{cor}\label{thm1.section}  Let $\pi$ be an irreducible orthogonal representation of $G$. Then $w(\pi)=1$.
\end{cor}

On the other hand, let $\rho$ be an irreducible symplectic representation of $G$. (By this we mean that $\rho$ is an irreducible representation on a complex vector space $V$ admitting a non-degenerate $G$-invariant \emph{antisymmetric} $B: V \times V \to \C$.) Its double $\pi=\rho \oplus \rho$ is orthogonal. Then
$\Theta_\pi(-\mathbbm{1})=-\Theta_\pi(\mathbbm{1})$, and so
\beq
\begin{split}
r_\pi &=\frac{1}{8}  \cdot 2 \Theta_\pi(\mathbbm{1}) \\
	&=\half \deg \rho. \\
	\end{split}
	\eeq
Therefore
\beq
w(\rho \oplus \rho)=(1+e)^{\half \deg \rho}.
\eeq

\subsection{$G=\SL(2,q)$, $q$ even}

Now let $q$ be even, say $q=2^r$.
 Let $N$ be the subgroup of upper unitriangular matrices in $G$. Being a Sylow 2-subgroup, $N$ detects the mod-2 cohomology of $G$.\\
Now $N\cong (\f_q,+)$ is an elementary abelian 2-group, so from Section \ref{cy}, we have
$$H^*(N)\cong H^*(C_2^r)\cong k[v_1,\hdots,v_r].$$
 Set $n_1=\begin{pmatrix}
1&&1\\
0&&1
\end{pmatrix}\in N.$

Every representation $\pi$ of $G$ is orthogonal by the main result in \cite{vinroot}.

\begin{thm} \label{sl2.qeven}
The total SWC of $\pi$ is
\beq
w(\pi) =\Big(\prod_{v\in H^1(N)}(1+v)\Big)^{s_\pi} \in H^*(N),
	\eeq
where $s_\pi=\frac1q(\Theta_\pi(1)-\Theta_\pi(n_1)).$
\end{thm}
\vspace{2mm}

The expansion of this product is well-known. We have
\begin{equation} \label{Dickson}
\prod_{v\in H^1(N)}(1+v)=1+\sum_{i=0}^{r-1}c_{r,i} \in H^*(N),\\
\end{equation}
with $c_{r,i}$ the famous Dickson invariants described in \cite{wilkerson1983primer}. Note that $\deg(c_{r,i})=2^r-2^i$.

\subsection{$G=\SL(3,q)$, $q$ odd}
Let $q$ be odd, and put $a_1=\diag(-1,-1,1)\in \SL(3,q)$. (Meaning, as usual, the diagonal matrix with these diagonal entries.)
\begin{thm}
  Let  $m_\pi=\frac18\big(\Theta_\pi(\mathbbm{1})-\Theta_\pi(a_1)\big).$ Then,
 \begin{equation*}
w(\pi)=
\begin{cases}
\big((1+t_1)(1+t_2)(1+t_1+t_2)\big)^{m_\pi}  \;\;, & q\equiv 1 \text{ (mod }4)\\
\big((1+v_1)(1+v_2)(1+v_1+v_2)\big)^{2m_\pi}\;,  & q\equiv 3 \text{ (mod }4).
\end{cases}
\end{equation*}
\end{thm}

\begin{remark} Again we encounter Dickson products as in \eqref{Dickson}, but where the vector space $H^1(N)$ is replaced with $\Span\{t_1,t_2\}$ or $H^1(A)$.
\end{remark}

\subsection{$G=\Sp(4,q)$, $q$ odd}
Let $q$ be odd, and $G=\Sp(4,q)$.  Put $g_1=\diag(1,-1,-1,1)\in \Sp(4,q).$
\begin{thm}
 The total SWC of an orthogonal representation $\pi$  of $G$ is
$$w(\pi)=((1+e_1)(1+e_2))^{r_\pi}(1+e_1+e_2)^{s_\pi},$$
where \begin{align*} r_\pi&=\frac{1}{16}\Big(\Theta_\pi(\mathbbm{1})-\Theta_\pi(-\mathbbm{1})\Big),\\
s_\pi&=\frac{1}{16}\Big(\Theta_\pi(\mathbbm{1})+\Theta_\pi(-\mathbbm{1})-2\Theta_\pi(g_1)\Big).
\end{align*}
\end{thm}

Again from Theorem 1 of \cite{gow1985real}, we deduce that when $\pi$ is irreducible orthogonal, we have
\beq
w(\pi)=(1+e_1+e_2)^{s_\pi},
\eeq
where 
\beq
s_\pi=\frac{1}{8} \left( \deg \pi-\Theta_\pi(g_1) \right).
\eeq

\subsection{$G=\Sp(6,q)$, $q$ odd}
Let \begin{align*}
g_1&=\text{diag}(1,1,-1,-1,1,1), \\
g_2&=\text{diag}(1,-1,-1,-1,-1,1)\in G.
\end{align*}
\begin{thm}
The total SWC of an orthogonal representation $\pi$ of $G$ is
$$w(\pi)=\Big(\prod_{1\leq i\leq 3}(1+e_i)\Big)^{r_\pi}\Big(\prod_{1\leq i<j\leq 3}(1+e_i+e_j)\Big)^{s_\pi}(1+e_1+e_2+e_3)^{t_\pi},$$
\end{thm}
where
\begin{align*}
r_\pi&=\frac{1}{32}\big(\Theta_\pi(\mathbbm{1})+\Theta_\pi(g_1)-\Theta_\pi(g_2)-\Theta_\pi(-\mathbbm{1})\big)\\
s_\pi&=\frac{1}{32}\big(\Theta_\pi(\mathbbm{1})-\Theta_\pi(g_1)-\Theta_\pi(g_2)+\Theta_\pi(-\mathbbm{1})\big), \text{ and}\\
t_\pi&=\frac{1}{32}\big(\Theta_\pi(\mathbbm{1})-3\Theta_\pi(g_1)+3\Theta_\pi(g_2)-\Theta_\pi(-\mathbbm{1})\big).
\end{align*}

Once more Theorem 1 of \cite{gow1985real} allows simplification: when $\pi$ is irreducible orthogonal, we have
\beq
w(\pi)=\Big(\prod_{1\leq i<j\leq 3}(1+e_i+e_j)\Big)^{s_\pi},
\eeq
where 
\beq
s_\pi=\frac{1}{16} \left( \deg \pi-\Theta_\pi(g_1) \right).
\eeq

\bibliographystyle{alpha}
\bibliography{mybib}
\vspace{10mm}
  
 \end{document}